# UNIQUENESS OF MAXIMAL ENTROPY MEASURE ON ESSENTIAL SPANNING FORESTS[1]

By Scott Sheffield

*Microsoft Research and University of California, Berkeley*

An *essential spanning forest* of an infinite graph $G$ is a spanning forest of $G$ in which all trees have infinitely many vertices. Let $G_n$ be an increasing sequence of finite connected subgraphs of $G$ for which $\bigcup G_n = G$. Pemantle's arguments imply that the uniform measures on spanning trees of $G_n$ converge weakly to an Aut($G$)-invariant measure $\mu_G$ on essential spanning forests of $G$. We show that if $G$ is a connected, amenable graph and $\Gamma \subset$ Aut($G$) acts quasitransitively on $G$, then $\mu_G$ is the unique $\Gamma$-invariant measure on essential spanning forests of $G$ for which the specific entropy is maximal. This result originated with Burton and Pemantle, who gave a short but incorrect proof in the case $\Gamma \cong \mathbb{Z}^d$. Lyons discovered the error and asked about the more general statement that we prove.

## 1. Introduction.

1.1. *Statement of result.* An *essential spanning forest* of an infinite graph $G$ is a spanning subgraph $F$ of $G$, each of whose components is a tree with infinitely many vertices. Given any subgraph $H$ of $G$, we write $F_H$ for the set of edges of $F$ contained in $H$. Let $\Omega$ be the set of essential spanning forests of $G$ and let $\mathcal{F}$ be the smallest $\sigma$-field in which the functions $F \to F_H$ are measurable.

Let $G_n$ be an increasing sequence of finite connected induced subgraphs of $G$ with $\bigcup G_n = G$. An Aut($G$)-invariant measure $\mu$ on $(\Omega, \mathcal{F})$ is Aut($G$)-*ergodic* if it is an extreme point of the set of Aut($G$)-invariant measures on $(\Omega, \mathcal{F})$. Results of [1, 8] imply that the uniform measures on spanning trees of $G_n$ converge weakly to an Aut($G$)-invariant and ergodic measure $\mu_G$ on $(\Omega, \mathcal{F})$.

Received July 2004; revised May 2005.

[1]Supported in part by NSF Grant DMS-04-03182.

We say $G$ is *amenable* if the $G_n$ above can be chosen so that

$$\lim_{n\to\infty}|\partial G_n|/|V(G_n)|=0,$$

where $V(G_n)$ is the vertex set of $G_n$ and $\partial G_n$ is the set of vertices in $G_n$ that are adjacent to a vertex outside of $G_n$. A subgroup $\Gamma\subset\mathrm{Aut}(G)$ *acts quasitransitively* on $G$ if each vertex of $G$ belongs to one of finitely many $\Gamma$ orbits. We say $G$ itself is *quasitransitive* if $\mathrm{Aut}(G)$ acts quasitransitively on $G$.

The *specific entropy* (also known as *entropy per site*) of $\mu$ is

$$-\lim_{n\to\infty}|V(G_n)|^{-1}\sum\mu(\{F_{G_n}=F_n\})\log\mu(\{F_{G_n}=F_n\}),$$

where the sum ranges over all spanning subgraphs $F_n$ of $G_n$ for which $\mu(\{F_{G_n}=F_n\})\neq 0$. This limit always exists if $G$ is amenable and $\mu$ is invariant under a quasitransitive action (see, e.g., [5, 7] for stronger results).

Let $\mathcal{E}_G$ be the set of probability measures on $(\Omega,\mathcal{F})$ that are invariant under some subgroup $\Gamma\subset\mathrm{Aut}(G)$ that acts quasitransitively on $G$ and that have maximal specific free entropy. Our main result is the following:

THEOREM 1.1. *If $G$ is connected, amenable and quasitransitive, then $\mathcal{E}_G=\{\mu_G\}$.*

1.2. *Historical overview.* As part of a long foundational paper on essential spanning forests published in 1993, Burton and Pemantle gave a short but incorrect proof of Theorem 1.1 in the case that $\Gamma\cong\mathbb{Z}^d$ and then used that theorem to prove statements about the dimer model on doubly periodic planar graphs [3]. In 2002, Lyons discovered and announced the error [6]. Lyons also extended part of the result of [3] to quasitransitive amenable graphs (Lemma 2.1 below) and questioned whether the version of Theorem 1.1 that we prove was true [6].

A common and natural strategy for proving results like Theorem 1.1 is to show first that each $\mu\in\mathcal{E}_G$ has a Gibbs property and second that this property characterizes $\mu$. The argument in [3] uses this strategy, but it relies on the incorrect claim that every $\mu\in\mathcal{E}_G$ satisfies the following property:

STRONG GIBBS PROPERTY. Fix any finite induced subgraph $H$ of $G$ and write $a\sim_O b$ if there is a path from $a$ to $b$ that consists of edges *outside* of $H$. Let $H'$ be the graph obtained from $H$ by identifying vertices equivalent under $\sim_O$. Let $\mu'$ be the measure on $(\Omega,\mathcal{F})$ obtained as follows: To sample from $\mu'$, first sample $F_{G\setminus H}$ from $\mu$ and then sample $F_H$ uniformly from the set of all spanning trees of $H'$. (We may view a spanning tree of $H'$ as a subgraph of $H$ because $H$ and $H'$ have the same edge sets.) Then $\mu'=\mu$. In other words, given $F_{G\setminus H}$—which determines the relation $\sim_O$ and the graph $H'$—the $\mu$ conditional measure on $F_H$ is the uniform spanning tree measure on $H'$.

This claim is clearly correct if $\mu = \mu_G$ and $G$ is a finite graph. To see a simple counterexample when $G$ is infinite, first recall that the number of *topological ends* of an infinite tree $T$ is the maximum number of disjoint semi-infinite paths in $T$ (which may be $\infty$). A *k-ended tree* is a tree with $k$ topological ends. If $G = \mathbb{Z}^d$ with $d > 4$, then $\mu_G \in \mathcal{E}_G$ and $\mu_G$-almost surely $F$ contains infinitely many trees, each of which has only one topological end [1, 8]. Thus, conditioned on $F_{G \setminus H}$, all configurations $F_H$ that contain paths joining distinct infinite trees of $F_{G \setminus H}$ have probability 0.

This example also shows, perhaps surprisingly, that $\mu \in \mathcal{E}_G$ does not imply that, conditioned on $F_{G \setminus H}$, all extensions of $F_{G \setminus H}$ to an element of $\Omega$ are equally likely. In other words, measures in $\mathcal{E}_G$ do not necessarily maximize entropy locally. Nonetheless, we claim that every $\mu \in \mathcal{E}_G$ does possess a Gibbs property of a different flavor:

WEAK GIBBS PROPERTY. For each $a$ and $b$ on the boundary of $H$, write $a \sim_I b$ if $a$ and $b$ are connected by a path contained *inside H* (a relationship that depends only on $F_H$). Then conditioned on this relationship and $F_{G \setminus H}$, all spanning forests $F_H$ of $H$ that give the same relationship (and for which each component of $F_H$ contains at least one point on the boundary of $H$) occur with equal probability.

If $\mu$ did not have this property, then we could obtain a different measure $\mu'$ from $\mu$ by first sampling a random collection $S$ of nonintersecting translates of $H$ (by elements of the group $\Gamma$) in a $\Gamma$-invariant way and then resampling $F_{H'}$ independently for each $H' \in S$ according to the conditional measure described above. It is not hard to see that $\mu'$ has higher specific entropy than $\mu$ and that it is still supported on essential spanning forests.

Unfortunately, the weak Gibbs property is not sufficient to characterize $\mu_G$. When $G = \mathbb{Z}^2$, for example, for each translation-invariant Gibbs measures on perfect matchings of $\mathbb{Z}^2$ there is a corresponding measure on essential spanning forests that has the weak Gibbs property [3]. The former measures have been completely classified and they include a continuous family of nonmaximal-entropy ergodic Gibbs measures [4, 9]. Significantly (see below), each of the corresponding nonmaximal-entropy measures on essential spanning forests almost surely contains infinitely many two-ended trees. The measure in which $F$ a.s. contains all horizontal edges of $\mathbb{Z}^2$ is a trivial example.

To prove Theorem 1.1, we will first show in Section 3.1 that if $\mu$ is $\Gamma$-invariant, has the weak Gibbs property and $\mu$-almost surely $F$ does not contain more than one two-ended tree, then $\mu = \mu_G$. We will complete the proof in Section 3.2 by arguing that if, with positive $\mu$ probability, $F$ contains more than one two-ended tree, then $\mu$ cannot have maximal specific entropy. Key elements of this proof include the weak Gibbs property, resamplings of $F$ on certain random extensions (denoted $\tilde{C}$ in Section 3.1) of finite subgraphs of $G$ and an entropy bound based on Wilson's algorithm.

We assume throughout the remainder of the paper that $G$ is amenable, connected and quasitransitive, $\Gamma$ is a quasitransitive subgroup of $\text{Aut}(G)$ and $G_n$ is an increasing sequence of finite connected induced subgraphs with $\bigcup G_n = G$ and $\lim |\partial G_n|/|V(G_n)| = 0$.

**2. Background results.** Before we begin our proof, we need to cite several background results. The following lemmas can be found in [3, 6, 8], [1, 3, 8] and [1, 2, 8], respectively.

LEMMA 2.1. *The measure $\mu_G$ is $\text{Aut}(G)$-invariant and ergodic, and has maximal specific entropy among quasi-invariant measures on the set of essential spanning forests of $G$. Moreover, this entropy is equal to*

$$-\lim_{n \to \infty} |V(G_n)|^{-1} \sum \mu_{G_n}(F_{G_n}) \log \mu_{G_n}(F_{G_n}),$$

*where $\mu_{G_n}$ is the uniform measure on all spanning forests $F_n$ of $G_n$ with the property that each component of $F_n$ contains at least one boundary vertex of $G_n$.*

LEMMA 2.2. *Let $C_n$ be any increasing sequence of finite subgraphs of $G$ whose union is $G$. For each $n$, let $H_n$ be an arbitrary subset of the boundary of $C_n$. Let $C'_n$ be the graph obtained from $C_n$ by identifying vertices in $H_n$. Then the uniform measures on spanning trees of $C'_n$ converge weakly to $\mu_G$. In particular, this holds for both wired boundary conditions $H_n = \partial C_n$ and free boundary conditions $H_n = \varnothing$.*

LEMMA 2.3. *If $G$ is amenable and $\mu$ is quasi-invariant, then $\mu$-almost surely all trees in $F$ contain at most two disjoint semi-infinite paths.*

We will also assume the reader is familiar with Wilson's algorithm for constructing uniform spanning trees of finite graphs by using repeated loop-erased random walks [10].

**3. Proof of the main result.**

3.1. *Consequences of the weak Gibbs property.*

LEMMA 3.1. *If $\mu$ has the weak Gibbs property and $\mu$-almost surely all trees in $F$ have only one topological end, then $\mu = \mu_G$.*

Proof. For a fixed finite induced subgraph $B$, we will show that $\mu$ and $\mu_G$ induce the same law on $F_B$. Consider a large finite set $C \subset V(G)$ that contains $B$. Then let $C_f$ be the set of vertices in $C$ that are starting points for infinite paths in $F$ that do not intersect $C$ after their first point. Then let $\tilde{C}$ be the union of $C_f$ and all vertices that lie on finite components of $F \backslash C_f$. In other words, $\tilde{C}$ is the set of vertices $v$ for which every infinite path in $F$ that contains $v$ includes an element of $C$.

Now, let $D$ be an even larger superset of $C$ that in particular contains all vertices that are neighbors of vertices in $C$. The weak Gibbs property implies that if we condition on the set $F_{G \backslash D}$ and the relationship $\sim_I$ defined using $D$, then all choices of $F_D$ that extend $F_{G \backslash D}$ to an essential spanning forest and preserve the relationship $\sim_I$ are equally likely. Now, if we further condition on the event $\tilde{C} \subset D$ and on a particular choice of $\tilde{C}$ and $C_f$, then all *spanning forests of $\tilde{C}$ rooted at $C_f$* (i.e., spanning trees of the graph induced by $\tilde{C}$ when it is modified by joining the vertices of $C_f$ into a single vertex) are equally likely to appear as the restriction of $F$ to $\tilde{C}$.

Since $D$ can be taken large enough so that it contains $\tilde{C}$ with probability arbitrarily close to 1, we may conclude that, in general, conditioned on $\tilde{C}$ and $C_f$, all spanning forests of $\tilde{C}$ rooted at $C_f$ are equally likely to appear as the restriction of $F$ to $\tilde{C}$. Since we can take $C$ to be arbitrarily large, the result follows from Lemma 2.2. □

Lemma 3.2. *If $\mu$ has the weak Gibbs property and $\mu$-almost surely $F$ consists of a single two-ended tree, then $\mu = \mu_G$.*

Proof. Define $B$ and $C$ as in the proof of Lemma 3.1. Given a sample $F$ from $\mu$, denote by $R$ the set of points that lie on the doubly infinite path (also called the *trunk*) of the two-ended tree. Then let $c_1$ and $c_2$ be the first and last vertices of $R$ that lie in $C$, and let $\tilde{C}$ be the set of all vertices that lie on the finite component of $F \backslash \{c_1, c_2\}$ that contains the trunk segment between $c_1$ and $c_2$. The proof is similar to that of Lemma 3.1, using the new definition of $\tilde{C}$ and noting that conditioned on $F_{G \backslash \tilde{C}}$, $c_1$ and $c_2$, all spanning trees of $\tilde{C}$ are equally likely to occur as the restriction of $F$ to $\tilde{C}$. The difference is that $\tilde{C}$ need not be a superset of $C$; however, we can choose a superset $C'$ of $C$ large enough so that the analogously defined $\tilde{C}'$ contains $C$ with probability arbitrarily close to 1. □

Lemma 3.3. *If $\mu$ has the weak Gibbs property and $\mu$-almost surely $F$ contains exactly one two-ended tree, then $\mu$ almost surely $F$ consists of a single tree and $\mu = \mu_G$.*

Proof. As in the previous proof, $R$ is the trunk of the two-ended tree. Clearly, each vertex in at least one of the $\Gamma$ orbits of $G$ has a positive probability of belonging to $R$. As in the previous lemmas, let $C$ be a large subset of $G$. Define $C_f$ to be the set of points in $C$ that are the initial points of infinite paths whose edges lie in the complement of $C$ and that belong to one of the single-ended trees of $F$. Let $\tilde{C}$ be the set of all vertices that lie on finite components of $F\backslash(C_f \cup \tilde{R})$. Conditioned on the trunk, $\tilde{C}$ and $C_f$, the weak Gibbs property implies that $F_{\tilde{C}}$ has the law of a uniform spanning tree on $\tilde{C}$ rooted at $\tilde{R} \cup C_f$ (i.e., vertices of that set are identified when choosing the tree).

Next we claim that if $R$ is chosen using $\mu$ as above, then a random walk started at any vertex of $G$ will eventually hit $R$ almost surely. Let $Q_R(v)$ be the probability, given $R$, that a random walk started at $v$ never hits $R$. Then $Q_R$ is harmonic away from $R$—that is, if $v \notin R$, then $Q_R(v)$ is the average value of $Q_R$ on the neighbors of $v$. If $v \in R$, then $Q_R(v) = 0$, which is at most the average value of $Q_R$ on the neighbors of $v$. Thus $Q(v) := \mathbb{E}_\mu Q_R(v)$ is subharmonic. Since $Q$ is constant on each $\Gamma$ orbit, it achieves its maximum, but if $Q$ achieves its maximum at $v$, it achieves a maximum at all of its neighbors and thus $Q$ is constant. Now, if $Q_R \neq 0$, then there must be a vertex $v$ incident to a vertex $w \in R$ for which $Q_R(v) \neq 0$, but then $Q_R(w)$ is strictly less than the average value at its neighbors: since $Q$ is harmonic, this happens with probability 0, and we conclude that $Q_R$ is $\mu$ a.s. identically 0.

It follows that if $C$ is a large enough superset of a fixed set $B$, then any random walk started at a point in $B$ will hit $R$ before it hits a point on the boundary of $C$ with probability arbitrarily close to 1. Letting $C$ get large (and choosing $C'$, as in the proof of the previous lemmma, large enough so that $\tilde{C}'$ contains $C$ with probability close to 1) and using Wilson's algorithm, we conclude that $\mu$-almost surely every point in $G$ belongs to the two-ended tree.  □

### 3.2. *Multiple two-ended trees.*

LEMMA 3.4. *If $\mu$ is quasi-invariant and with positive $\mu$ probability $F$ contains more than one two-ended tree, then the specific entropy of $\mu$ is strictly less than the specific entropy of $\mu_G$.*

PROOF.    Let $k$ be the smallest integer such that for some $v \in V(G)$, there is a positive $\mu$ probability $\delta$ that $v$ lies on the trunk $R_1$ of a two-ended tree $T_1$ of $F$ and is distance $k$ from the trunk $R_2$ of another two-ended tree of $F$. We call a vertex with this property a *near intersection* of the ordered pair $(R_1, R_2)$. Let $\Theta$ be the $\Gamma$ orbit of a vertex with this property. Every $v \in \Theta$ is a near intersection with probability $\delta$.

Flip a fair coin independently to determine an orientation for each of the trunks. Fix a large connected subset $C$ of $G$. Let $C_f$ be the set containing the last element of each component of the intersection of $C$ with a trunk and let $C_b$ be the set of all of the first elements of these trunk segments. Let $\overline{C}_f$ be the union of $C_f$ and one vertex of $\partial C$ from each tree of $F_C$ that does not contain a segment of a trunk. We may then think of $F_C$ as a spanning forest of the graph induced by $C$ rooted at the set $\overline{C}_f$.

Let $\nu$ be the uniform measure on *all* spanning forests of $C$ rooted at $\overline{C}_f$. Denote by $C^k$ the set of vertices in $C \cap \Theta$ of distance at least $k$ from $\partial C$. Let $A = A(C, C_b, \overline{C}_f, m)$ be the event that the paths from $C_b$ to $\overline{C}_f$ are disjoint paths that end at the $C_f$ and have at least $m$ near intersections in $C^k$. We will now give an upper bound on $\nu(A)$ (which is zero if either $C_b$ or $\overline{C}_f$ is empty).

We can sample from $\nu$ using Wilson's algorithm, beginning by running loop-erased random walks starting from each of the points in $C_b$ to generate a set of paths from the points in $C_b$ to the set $\overline{C}_f$ (which may or may not join up before hitting $\overline{C}_f$). Order the points in $C_b$ and let $P_1, P_2, \ldots$ be the paths beginning at those points. For any $r, s \geq 1$, Wilson's algorithm implies that conditioned on $P_i$ with $i < r$ and on the first $s$ points $P_r$, the $\nu$ distribution of the next step of $P_r$ is that of the first step of a random walk in $C$ beginning at $P_r(s)$ and conditioned not to return to $P_r(1), \ldots, P_r(s)$ before hitting either $\overline{C}_f$ or some $P_i$ with $i < r$.

For each $r > 1$, we define the first *fresh near collision point* (FNCP) of $P_r$ to be the first point in $P_r$ that lies in $C^k$ and is distance $k$ or less from a $P_i$ with $i < r$. The $j$th FNCP is the first point in $P_r$ that lies in $C^k$, is distance $k$ or less from a $P_i$ with $i < r$ and is distance at least $k$ from the $(j-1)$st FNCP in $P_r$. If we condition on the $P_1, P_2, \ldots, P_{r-1}$ and on the path $P_r$ up to an FNCP, then there is some $\varepsilon$ (independent of details of the paths $P_i$) such that with $\nu$ probability at least $\varepsilon$, after at most $k$ more steps, the path $P_r$ collides with one of the other $P_i$. Let $K$ be the total number of vertices of $G$ within distance $k$ of a vertex $v \in \Theta$. Since on the event $A$, we encounter at least $m/K$ FNCP's (as every near intersection lies within $k$ units of an FNCP) and the collision described above fails to occur after each of them, we have $\nu(A) \leq (1 - \varepsilon)^{m/K}$.

Let $B = B(n, m) \in \mathcal{F}$ be the event that when $C = G_n$, $F_C \in A(C, C_b, \overline{C}_f, m)$ for *some* choice of $C_b$ and $\overline{C}_f$. Summing over all the choices of $\overline{C}_f$ and $C_b$ (the number of which is only exponential in $|\partial G_n|$), we see that if $m$ grows linearly in $|V(G_n)|$, then $\mu_{G_n}(B(n, m))$ (where $\mu_{G_n}$ is defined as in Lemma 2.1) decays exponentially in $|V(G_n)|$. [Note that since $\nu$ is the uniform measure on a subset of the support of $\mu_{G_n}$, any $X$ in the support of $\nu$ has $\mu_{G_n}(X) \leq \nu(X)$.]

Because the expected number of near collisions is linear in $|V(G_n)|$, there exist constants $\varepsilon_0$ and $\delta_0$ such that for large enough $n$, there are at least $\delta_0 |V(G_n)|$ near intersections in $G_n^k$ with $\mu$ probability at least $\varepsilon_0$. However, the $\mu_{G_n}$ probability that this occurs decays exponentially in $|V(G_n)|$. From this, it is not hard to see that the specific entropy of the restriction of $\mu$ to $G_n$ [i.e., $-|V(G_n)|^{-1} \sum \mu(F_{G_n}) \log \mu(F_{G_n})$] is less than the specific entropy of $\mu_{G_n}$ [i.e., $|V(G_n)|^{-1} \log N$, where $N$ is the size of the support of $\mu_{G_n}$] by a constant independent of $n$. By Lemma 2.1, the specific entropy of $\mu_{G_n}$ converges to that of $\mu_G$, so the specific entropy of $\mu$ must be strictly less than that of $\mu_G$. □

**Acknowledgments.** We thank Russell Lyons for suggesting the problem, for helpful conversations, and for reviewing early drafts of the paper. We also thank Oded Schramm and David Wilson for helpful conversations.

MATHEMATICS DEPARTMENT
UNIVERSITY OF CALIFORNIA
BERKELEY, CALIFORNIA 94720-3860
USA
E-MAIL: sheff@math.berkeley.edu